\documentclass[12pt]{article}
\usepackage{amssymb}

\def\appendix{\par}  

\def\cc{{\mathfrak c}}
\def\ff{{\mathcal F}}
\def\ii{{\mathcal I}}
\def\uu{{\mathcal U}}
\def\vv{{\mathcal V}}
\def\bb{{\mathcal B}}
\def\pp{{\mathbb P}}
\def\infsub{[\omega]^\omega}

\def\bb{{\mathfrak b}}
\def\term{{\rm term}}
\def\proof{\par\noindent Proof\par\noindent}

\def\forces{{| \kern -2pt \vdash}}

\def\qed{\par\noindent QED\par}
\def\rmand{\mbox{ and }}
\def\rmor{\mbox{ or }}
\def\rmiff{\mbox{ iff }}
\def\sig{{\bf\Sigma}^1_1}
\def\pioneone{{\bf\Pi}^1_1}
\def\pionetwo{{\Pi}^1_2}

\newtheorem{theorem}{Theorem}
\newtheorem{lemma}[theorem]{Lemma}

\newtheorem{question}[theorem]{Question}

\newtheorem{prop}[theorem]{Proposition}

\begin{document}

\begin{center}
{\large Ultrafilters with property (s)}
\end{center}

\begin{flushright}
Arnold W. Miller\footnote{
Thanks to the Fields Institute, Toronto for their support
during the time these results were proved and to Juris Steprans
for helpful conversations and thanks to Boise State University  
for support during the time this paper
was written. 
}
\end{flushright}

\begin{center}
Abstract
\end{center}

\begin{quote}

A set $X\subseteq 2^\omega$ has property (s) (Marczewski (Szpilrajn)) iff
for every perfect set $P\subseteq 2^\omega$ there exists a perfect set
$Q\subseteq P$ such that $Q\subseteq X$ or $Q\cap X=\emptyset$.
Suppose $\uu$ is a nonprincipal ultrafilter on $\omega$.
It is not difficult to see that if $\uu$ is preserved by Sacks forcing,
i.e., it generates an ultrafilter in the generic extension after forcing
with the partial order of perfect sets, then $\uu$ has property (s) in
the ground model.  It is known that selective ultrafilters or even
P-points are preserved by Sacks forcing.  On the other hand (answering
a question raised by Hrusak) we show that assuming CH (or more generally
MA) there exists an ultrafilter $\uu$ with property (s) such that 
$\uu$ does not generate an ultrafilter in any extension which adds a new
subset of $\omega$.

\end{quote}

It is a well known classical result due to Sierpinski (see \cite{bj}) 
that a nonprincipal ultrafilter $\uu$ on $\omega$ when considered as a
subset of $P(\omega)=2^\omega$ cannot have the property of Baire or
be Lebesgue measurable.  Here we identify 
$2^\omega$ and $P(\omega)$ by identifying a subset of $\omega$ with
its characteristic function.
Another very weak regularity property
is property (s) of Marczewski (see Miller \cite{survey}). A set
of reals $X\subseteq 2^\omega$ has property (s) iff for every perfect
set $P$ there exists a subperfect set $Q\subseteq P$ such that
either $Q\subseteq X$ or $Q\cap X=\emptyset$.  Here by perfect we mean
homeomorphic to $2^\omega$.

It is natural to ask:

\bigskip\noindent
{\bf Question.} (Steprans) Can a nonprincipal ultrafilter
$\uu$ have property (s)?

\bigskip
If $\uu$ is an ultrafilter in a model of set theory $V$ and 
$W\supseteq V$ is another model of set theory then we say
$\uu$ generates an ultrafilter in $W$ if for every
$z\in P(\omega)\cap W$ there exists $x\in\uu$ with
$x\subseteq z$ or $x\cap z=\emptyset$.  This means that the
filter generated by $\uu$ (i.e. closing under supersets) 
is an ultrafilter in $W$.

We begin with the following result:

\begin{theorem}\label{one}
For $\uu$ a nonprincipal ultrafilter on $\omega$ in $V$ the following are
equivalent:
\begin{enumerate}
 \item For some Sack's generic real $x$ over $V$
 $$V[x]\models \uu \mbox{ generates an ultrafilter. }$$
 \item In $V$, for every perfect set $P\subseteq P(\omega)$ there 
 exists a perfect set 
 $Q\subseteq P$
 and a $z\in \uu$ such that either $\forall x\in Q\;\; z\subseteq x$
 or $\forall x\in Q\;\; z \cap x=\emptyset$.
 \item For some extension $W\supseteq V$ with a new subset of $\omega$
 $$W\models \uu \mbox{ generates an ultrafilter. }$$
\end{enumerate}
\end{theorem}

\proof

To see that $(3)\to (2)$, let $P$ be any perfect set coded in $V$.   Since $W$
contains a new subset of $\omega$ there exists $x\in (P\cap W)\setminus V$.  

Since $\uu$ generates an ultrafilter in $W$ there  exists $z\in \uu$ so that
either $z\subseteq x$ or $z\cap x=\emptyset$.  Suppose the first happens. In
$V$ consider the set 
 $$Q=\{y\in P: z\subseteq y\}$$ 
Note that the new real $x$ is in the closed set $Q$.  It follows that $Q$ must
be an uncountable closed set and so it contains a perfect subset. The other
case is exactly the same.  

One way to see that $Q$ must be  uncountable is to note that if (in $V$)
$Q=\{x_n:n<\omega\}$, then the $\pioneone$ sentence  
$$\forall x\in 2^\omega (x\in Q\rmiff \exists n<\omega \;x=x_n)$$ 
would be true in $V$ and since $\pioneone$
sentences are absolute (Mostowski absoluteness, see \cite{kechris}) true in
$W$.  Another way to prove it is to do the standard derivative Cantor argument
to the closed set $Q$ removing isolated points and
iterating thru the transfinite and noting that each real removed is in $V$, 
while the new
real is never removed, and hence the kernel of $Q$ is perfect. See Solovay
\cite{solovay}, for a similar proof of Mansfield's theorem that a (lightface)
$\pionetwo$ set with a nonconstructible element contains a perfect set.

Now we see that $(2)\to (1)$.
A basic property of Sack's forcing is that for every $y\in 2^\omega \cap M[x]$
in a Sacks extension is either in $M$ or is itself Sacks generic over $M$ (see
Sacks \cite{sacks}).  Hence we need only show that if $y\subseteq \omega$ is
Sacks generic over $M$, then there exists $z\in \uu$ with $z\subseteq y$
or $z\cap y=\emptyset$.  Recall also that the Sacks
real $y$ satisfies that the generic filter $G$ is
exactly the set of all perfect sets $Q$ coded in $V$ with $y\in Q$. 
 
Condition (2) says that the set of such $Q$ are dense and hence there exists
$Q$ in the generic filter determined by $y$ and $z\in \uu$ such that either 
$z\subseteq u$ for every $u\in Q$  or  either $z\cap u=\emptyset$ for every
$u\in Q$.  But this means that either $z\subseteq y$ or $z\cap y=\emptyset$.

$(1)\to (3)$ is obvious.  

\qed

\bigskip
Remark.
The above proof also shows that if an ultrafilter is preserved
in one Sacks extension, then it is preserved in all Sacks extensions.

\bigskip
Remark.  In Baumgartner and Laver \cite{bl} it is shown that selective
ultrafilters are preserved by Sacks forcing.  In Miller \cite{super}
it is shown that $P$-points are preserved by superperfect set forcing
(and hence by Sacks forcing also). 

\bigskip
We say that an ultrafilter $\uu$ is preserved by Sacks forcing iff
for some (equivalently all) Sacks generic reals $x$ that 
$\uu$ generates an ultrafilter in $V[x]$.  Recall that $\uu\times \vv$
is the ultrafilter on $\omega\times\omega$ defined by
$$A\in \uu\times \vv \rmiff \{n:\{m:(n,m)\in  A\}\in\uu\}\in\vv$$
If $\uu$ and $\vv$ are nonprinciple ultrafilters, then $\uu\times\vv$
is not a P-point.  Also recall that $\uu\leq_{RK}\vv$ (Rudin-Keisler)
iff there exists $f\in\omega^\omega$ such that for every $X\subseteq \omega$
$$X\in\uu \rmiff f^{-1}(X)\in \vv$$

\bigskip
\begin{prop}\label{two}
If $\uu$ and $\vv$ are preserved by Sacks forcing, then so is $\uu\times\vv$.
If $\uu\leq_{RK}\vv$ and $\vv$ is preserved by Sacks forcing, then so
is $\uu$.
\end{prop}
\proof
Suppose $A\subseteq \omega\times\omega$ and $A\in V[x]$.  For
each $n<\omega$ let $A_n=\{m:(n,m)\in A\}$.  Since $\uu$ is preserved
there exists $B_n\in\uu$ with $B_n\subseteq A_n$ or $B_n\cap A_n=\emptyset$.
By the preservation of $\vv$ there exists $C\in\vv$ such that either
$B_n\subseteq A_n$ for all $n\in C$ or
$B_n\cap A_n=\emptyset$ for all $n\in C$.  
By the Sacks property there exists 
$(\bb_n\in [\uu]^{2^n}:n<\omega)\in V$ such that $B_n\in\bb_n$ 
for every $n$. Let $B_n^0=\cap \bb_n$.  Then
$$\cup_{n\in C}\{n\}\times B_n^0\subseteq A \rmor
\cup_{n\in C}\{n\}\times B_n^0\cap A=\emptyset
$$ 

Suppose $\uu\leq_{RK}\vv$ via $f$.  If $A\subseteq \omega$, then
since $\vv$ is preserved, there exists $B\in\vv$ such that
either $B\subseteq f^{-1}(A)$ or $B\subseteq f^{-1}(\overline{A})$.
but then $f(B)\subseteq A$ or $f(B)\subseteq \overline{A}$
and since $f(B)\in\uu$ we are done.
\qed
 
\bigskip Remark. The Rudin-Keisler result is generally true, but the
product result depends on the bounding property.  For example, if $\uu$ is a
P-point, then $\uu$ is preserved in the superperfect extension, but
$\uu\times\uu$ is not.

\bigskip

It is clear that property (2) of Theorem \ref{one} 
implies that any ultrafilter which is
preserved by Sacks forcing has property (s).  But what about the converse?
The main result of this paper is that the reverse implication is false.
This answers a question raised by Hrusak.

\begin{theorem}
Suppose the CH is true or even just that the real line cannot be
covered by fewer than continuum many meager sets.
Then there exists an ultrafilter $\uu$ on
$\omega$ which has property (s) but is not preserved by Sacks forcing.
\end{theorem}

\proof

We give the proof in the case of the continuum hypothesis and indicate
how to do it under the more general hypothesis.

Let $\ii\subseteq\infsub$ be an independent perfect family. 
Independent means that for every $m,n$ and distinct 
$x_1,\ldots,x_m,y_1,\ldots,y_n\in\ii$ the set
$$x_1\cap\ldots\cap x_m\cap \overline{y}_1
\cap \ldots\cap \overline{y}_n\mbox{ is infinite.}$$
where $\overline{y}$ means the complement of $y$ in $\omega$.
We claim that the following family 
$$\ii\cup\{\overline{z}:\;\exists^\infty \; x\in \ii\;\; z\subseteq^*x\}$$
has the finite intersection property.  ($\exists^\infty$ means
there exists infinitely many).  To see this suppose that
$x_1,\ldots,x_m\in \ii$ and $z_1,\ldots, z_n$ and 
$\exists^\infty \; x\in \ii\;\; z_i\subseteq^*x$ for each $i$.  Then
we can choose $y_i\in \ii$ distinct from each other and the $x's$ so
that each $z_i\subseteq^* y_i$.  But since $\overline{y}_i\subseteq
\overline{z}_i$ we have that

$$x_1\cap\ldots\cap x_m\cap \overline{y}_1 \cap \ldots\cap \overline{y}_n
\subseteq^* x_1\cap\ldots\cap x_m\cap \overline{z}_1 \cap \ldots\cap
\overline{z}_n $$ 
By independence the set on the left is infinite and hence so
is the set on the right. 
Thus this family has the finite intersection property.

Now let
$\ff_0$ be the filter generated by 
$\ii\cup\{\overline{z}:\exists^\infty \; x\in \ii\;\; z\subseteq^*x\}$.

Note that if $\uu\supseteq\ff_0$ is any ultrafilter then it cannot be
preserved by Sacks forcing.  This is because $\ii$ is a perfect subset
of $\uu$, however there is no $z\in \uu$ with $z\subseteq x$ for
all $x\in \ii$ or even infinitely many $x\in\ii$
or else $\overline{z}\in\ff_0\subseteq\uu$, hence Theorem \ref{one} (2)
fails.

Note that since $\ii$ was perfect the filter $\ff_0$ is a $\sig$
subset of $P(\omega)$. 

\begin{lemma}
Suppose that $P\subseteq\infsub$ is a perfect set and 
$\ff$ is a $\sig$ filter extending the cofinite filter 
on $\omega$.  Then there exists a perfect $Q\subseteq P$ such that
either
\begin{enumerate}
 \item $\ff\cup Q$ has the finite intersection property or
 \item there exists $z\subseteq \omega$ so that $\ff\cup\{z\}$ has
 the finite intersection property and for every $x\in Q$ we have that
 there exists $y\in\ff$ with $x\cap y\cap z=\emptyset$. 
\end{enumerate}
\end{lemma}
\proof
The strategy is try to do a fusion argument to get case (1).   If it
every fails, then stop and get case (2).

\bigskip {\bf Claim}. Suppose $(Q_i:i<n)$ are disjoint perfect subsets of
$\infsub$  Then either there exists $(Q_i'\subseteq Q_i:i<n)$ perfect so that
for every $(x_i\in Q_i':i<n)$ and $y\in \ff$ we have that  
 $$|y\cap x_0\cap x_1\cap \ldots \cap x_{n-1}|=\omega$$ 
or there exists $z\subseteq \omega$, $k<n$,
and $Q\subseteq Q_k$ perfect so that $\ff\cup\{z\}$ has the finite intersection
property and for every $x\in Q$ we have that there exists $y\in\ff$ with $x\cap
y\cap z=\emptyset$.
\proof
Consider 
$$A_k=\{(x_i\in\infsub:i<k): \exists y\in\ff\;\; |y\cap x_0\cap x_1
\cap\ldots x_{k-1}|<\omega\}$$
Since $\ff$ is $\sig$ it is easy to see that each $A_k$ is a 
$\sig$ set and hence has the property of Baire relative to
the product $\prod_{i<k}Q_i$.  
By Mycielski \cite{mycielski} (see also Blass \cite{blass})
there exists perfect sets 
$(Q_i^*\subseteq Q_i:i<n)$ so that for every $k\leq n$ either 
$$\prod_{i<k}Q_i^*\cap A_k=\emptyset\rmor\prod_{i<k}Q_i^*\subseteq A_k.$$  
If the first case happens for $k=n$, then we let $Q_i'=Q_i^*$ and 
the claim is proved. If the second case happens choose $k$ minimal for which
it happens.  This means we have that 
\begin{enumerate}
 \item for all $(x_i:i<k-1)\in \prod_{i<k-1} Q^*_i$ and 
 $y\in\ff$ we have $$|y\cap x_0\cap x_1\cap \ldots \cap x_{k-2}|=\omega$$
 and
 \item for all $(x_i:i<k)\in \prod_{i<k} Q^*_i$ there exists
 $y\in\ff$ such that 
 $$|y\cap x_0\cap x_1\cap \ldots \cap x_{k-1}|=\emptyset$$
\end{enumerate}  
In this case let $(x_i:i<k-1)\in \prod_{i<k-1} Q^*_i$ be arbitrary and
put 
$$z=x_0\cap x_1\cap \ldots \cap x_{k-2}\rmand Q=Q^*_k\subseteq Q_k.$$  
This proves the Claim.
\qed
It is now an easy fusion argument to finish proving the Lemma from the
Claim.
\qed

Now we construct our ultrafilter proving the theorem under the assumption
of CH.  
We let $(P_\alpha:\alpha<\omega)$ list all perfect subsets of $2^\omega$.
We construct an increasing sequence $\ff_\alpha$ for $\alpha<\omega_1$ of
$\sig$ filters as follows. 

Let $\ff_0$ the filter generated by 
$$\ii\cup\{\overline{z}:\exists^\infty \; x\in \ii\;\; z\subseteq^*x\}$$
At limit ordinals $\alpha$ we let $\ff_\alpha$ be the union 
$\cup_{\beta<\alpha}\ff_beta$ and note that it is a $\sig$ filter.
At successor stages $\alpha+1$ we apply the Lemma to $P_\alpha$
and $\ff_\alpha$.  In the first case we find a perfect set 
$Q\subseteq P_\alpha$ such that $\ff_\alpha\cup Q$ has the finite
intersection property.  In this case we let $\ff_{\alpha+1}$
be the filter generated by $\ff_\alpha\cup Q$ and note that is
$\sig$.  In the second case we find a perfect set $Q\subseteq P_\alpha$
and $z\subseteq\omega$ so that $\ff_\alpha\cup\{z\}$ has 
the finite intersection
property and for every $x\in Q$ we have that
there exists $y\in\ff_\alpha$ with $x\cap y\cap z=\emptyset$.  Here we let
$\ff_{\alpha+1}$ be the filter generated by $\ff_\alpha\cup\{z\}$ and
note that for every ultrafilter $\uu\supseteq \ff_{\alpha+1}$ that
$\uu\cap Q=\emptyset$, because we have put 
$\{\overline{x}:x\in Q\}\subseteq\ff_{\alpha+1}$.  this ends the proof
under CH.

Now we see how do this construction under the weaker hypothesis that  the real
line cannot be covered by fewer than continuum many meager sets.
We construct an increasing sequence $(\ff_\alpha:\alpha<\cc)$ of filters
such that each $\ff_\alpha$ is the union of $\leq |\alpha|$ $\;\;\;\sig$
sets.  In order to prove the corresponding Claim and Lemma we note that
the following is true.  

\bigskip{\bf Claim} Suppose the real line cannot be covered by 
$\kappa$ many meager sets, $(Q_k:k<n)$ are perfect, and for
each $k\leq n$ we have
$A_k\subseteq\prod_{i<k}\infsub$ which is the union of $\leq\kappa$
many $\sig$ sets.  Then there exists  
$(Q_i^*\subseteq Q_i:i<n)$ perfect so that for every $k\leq n$ either
$$\prod_{i<k}Q_i^*\subseteq A_k\rmor \prod_{i<k}Q_i^*\cap A_k=\emptyset.$$

\proof
Construct $(Q_i^j:i<n)$ perfect by induction so that
\begin{enumerate}
 \item $Q_i^0=Q_i$ all $i<n$, 
 \item $Q_i^{j+1}\subseteq Q_i^{j}$ all $i<n$,
 \item $\prod_{i<j}Q_i^j\subseteq A_j$ or $\prod_{i<j}Q_i^j\cap A_j=\emptyset$.
\end{enumerate}
Given $(Q_i^j:i<n)$ and $A_{j+1}$ the union of $\kappa$ many
$\sig$ sets, say $\cup \{B_\alpha:\alpha<\kappa\}$ there are two
cases. 

Case 1. For some $\alpha<\kappa$ the set
$B_\alpha\cap  \prod_{i<j+1}Q_i^j$ is not meager.

In this case it must be comeager in some relative interval
$\prod_{i<j+1}Q_i^j\cap \prod_{i<j+1}[s_i]$. And now we can find
$Q_i^{j+1}\subseteq  Q_i^j\cap [s_i]$ such that  
$$\prod_{i<j+1}Q_i^{j+1}\subseteq B_\alpha\subseteq A_{j+1}$$

Case 2. Each $B_\alpha$ is meager in $\prod_{i<j+1}Q_i^j$.  

In this case we use the covering of category hypothesis in the form
of Martin's axiom for countable posets.  For $p\subseteq \omega^{<\omega}$
a finite subtree, define $s$ a terminal node of $p$ (
$s\in\term(p)$) iff $s\in p$ and for every $t\in p$ if
$s\subseteq t$ then $s=t$. 
For $p,q$ finite subtrees
of $2^{<\omega}$ we define $p\supseteq_{e}q$ (end extension) iff
$p\supseteq q$ and every new node of $p$ extends a terminal node of
$q$.  Define $T_i=\{s\in 2^{<\omega}: [s]\cap Q_i\not=\emptyset\}$.

Consider the partial
order $\pp$ consisting of finite approximations to products of
perfect trees below the $Q_i$:
$$\pp=\{(p_i:i<n): p_i \mbox{ is a finite subtree of } T_i, i<n\}$$
and $p\leq q$ iff $p_i\supseteq_e q_i$ all $i<n$.  
Our assumption about the covering of the real by meager sets is equivalent
to MA$_\kappa($ctble), i.e., for every countable poset $\pp$ and any
family $(D_\alpha:\alpha<\kappa)$ of dense subsets of $\pp$ there exists
a $\pp$-filter $G$ such that $G\cap D_\alpha\not=\emptyset$ for all
$\alpha<\kappa$.  Note that for any $k<n$ and $C\subseteq \prod_{i<k}Q_i$
which is nowhere dense the set
$$D_C=\{p\in\pp\;:\; \forall i<k\;\;\forall (s_i\in\term(p_i):
i<k)\;\; C\cap \prod_{i<k}[s_i]=\emptyset\}$$
is dense in $\pp$ similarly, for any $m<\omega$ the following sets 
are dense:
$$D_{m}=\{p\in\pp\;:\; \forall i<n\forall s\in\term(p_i)\;\;|s|>m\}$$
$$D_m^*=\{p\in\pp\;:\; \forall i<n\forall s\in p_i\;\;|s|=m
\to \exists t \; t\supseteq s \rmand t0,t1\in p_i\}$$

So a sufficiently generic filter produces a sequence 
$(T_i'\subseteq T_i:i<n)$ of perfect subtrees such
that letting $Q_i'=[T_i']$ we with the property that
$\prod_{i<n}Q_i'\cap B_\alpha=\emptyset$ for all $\alpha<\kappa$.
\qed

\begin{question}
Can we prove in ZFC that there exists a nonprinciple ultrafilter with
property (s)? 
\end{question}

\bigskip
Remark. It is easy to construct a nonprinciple ultrafilter which fails
to have property (s).  Start with a perfect independent family 
$I\subseteq P(\omega)$. Choose 
$$\{X_\alpha:\alpha<\cc\}\cup \{Y_\alpha:\alpha<\cc\}\subseteq I$$
distinct so that for every perfect $Q\subseteq I$ there exists
$\alpha$ with $X_\alpha$ and $Y_\alpha$ both in $Q$.  Then any
ultrafilter 
$$\uu\supseteq
\{X_\alpha:\alpha<\cc\}\cup \{\overline{Y}_\alpha:\alpha<\cc\}$$
will fail to have property (s).

\bigskip

\begin{question}
Can we prove in ZFC that there exists a nonprinciple ultrafilter with
which is preserved by Sacks forcing? 
\end{question}

Note Shelah,  see  \cite{bj}, has shown it is consistent that there are no
nonprinciple P-points.  See Brendle \cite{brend} for a plethora of
ultrafilters weaker than P-points such as Baumgartner's nowhere dense
ultrafilters.

\begin{question}
Suppose $\uu$ and $\vv$ are nonprinciple ultrafilters in $V$ which generate
ultrafilters  $\uu^*$ and $\vv^*$ in $W\supseteq V$.  
If $\uu^*\leq_{RK}\vv^*$ holds in $W$, must $\uu\leq_{RK}\vv$ 
be true in $V$? 
\end{question}

\begin{question}
Suppose $\uu\in V$ is generates an ultrafilter
in $V[x]$ for some (equivalently all)  Sacks reals $x$ over $V$.  Suppose 
$x$ is
a Sacks real over $V$ and $y$ is a Sacks real  over $V[x]$.  Must
$\uu$ generate an ultrafilter in $V[x,y]$?
\end{question}

\begin{question}
Suppose $V\subseteq W$ and 
$$V\models \uu \mbox{ has property (s)}$$
does 
$$W\models \{A\subseteq \omega:\exists B\in\uu \;B\subseteq A\}=\uu^*
\mbox{ has property (s) ?}$$
\end{question}

\begin{question}
Is Proposition \ref{two} true for property (s) in place of ``preserved
by Sacks forcing''?
\end{question}

\begin{flushleft}
Arnold W. Miller \\
miller@math.wisc.edu \\
http://www.math.wisc.edu/$\sim$miller\\
University of Wisconsin-Madison \\
Department of Mathematics, Van Vleck Hall \\
480 Lincoln Drive \\
Madison, Wisconsin 53706-1388 \\
\end{flushleft}

\end{document}